\documentclass{article}

\usepackage{amsmath}
\usepackage{amssymb}
\usepackage[applemac]{inputenc}

\newcommand{\Weyl}{\mathrm{Weyl}}
\newcommand{\Sq}{\mathrm{Sq}}
\newcommand{\Inv}{\mathrm{Inv}}
\newcommand{\cG}{\mathcal{G}}

\newcommand{\Hom}{\mathrm{Hom}}
\newcommand{\Aut}{\mathrm{Aut}}
\newcommand{\ch}{\mathrm{char}}

\newcommand{\Gal}{\mathrm{Gal}}
\newcommand{\n}{\noindent}

\newcommand{\GL}{\mathrm{GL}}

\newcommand{\Or}{\mathrm{O}}

\newcommand{\Z}{\mathbf Z}

\newcommand{\F}{\mathbf F}
\newcommand{\R}{\mathbf R}

\newcommand{\w}{w^{\mathrm gal}}
\newcommand{\e}{\mathrm{e}}

\begin{document}

\title{Lines on cubic surfaces,  Witt invariants and Stiefel-Whitney classes}
\author{Eva Bayer-Fluckiger and Jean-Pierre Serre}
\date{\today}
\maketitle


\medskip


{\bf Introduction} 

\medskip
  This paper is a contribution to the classical topic of the 27 lines of a cubic surface,
  and the Weyl group $G$ of $\sf{E}_6$.
  
  We show how the theory of ``cohomological invariants" and ``Witt invariants" can be used to clarify the structure of some of the quadratic forms
  and Stiefel-Whitney classes associated with the surface.
  
  \smallskip

The first two sections recall well-known facts about several lattices associated with $G$. These lattices are $G$-stable, hence give rise to {\it $G$-quadratic forms} over any ground field $k$ of characteristic $\neq 2$. Among these,
the quadratic form of rank 7 has a geometric interpretation in terms of sheaf cohomology (Levine-Raksit [LR 18], cf. $\S$7.3), and the rank $27$ one is the trace form
of the \'etale algebra of the $27$ lines. The relations between these forms are discussed in $\S$$\S$3-6.
The main tool of the proofs is the ``splitting principle" for Witt invariants, as in [Se 03]. 

\medskip

An appendix gives a formula for
  the Stiefel-Whitney classes of any linear representation of $G$, in terms of four basic invariants $w_1,w_2,w_3,w_4$. Here also, the splitting principle
  plays an essential role.

\medskip
{\bf \S1. Lattices associated with} $\Weyl(\sf{E}_6)$. 

\medskip

  Let $R$ be a root system of type $\sf{E}_6$, with basis $\{\alpha_1,...,\alpha_6\}$ numbered as in Bourbaki [Bo 68], \S4.12:

  \medskip
  
\hspace{25mm}$\alpha_{_1} \hspace{6mm}\alpha_{_3}\hspace{5mm}\alpha_{_4}\hspace{6mm} \alpha_{_5}\hspace{5mm} \alpha_{_6}$

\hspace{25mm}$\circ$-------$\circ$------$\circ$-------$\circ$------$\circ$

\hspace{43mm} $\mid$

\hspace{43mm} $\circ \ \alpha_{_2}$

\hspace{43mm} $\vdots $

\hspace{43mm} $\circ  \ \alpha_{_0}$

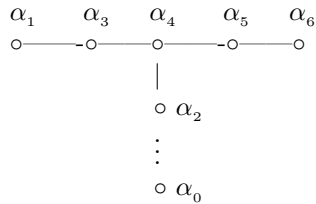
\smallskip

\n where $\alpha_{_0}$ is the opposite of the highest root $\tilde{\alpha}= \alpha_1+2\alpha_2+2\alpha_3+3\alpha_4+2\alpha_5+\alpha_6$. Let
  $\{\omega_1,...,\omega_6\}$ be the corresponding fundamental weights. Let $Q$ be the root lattice, which is generated by the $\alpha_i$; its discriminant is 3.
  Its dual lattice $P$ is the weight lattice, with basis the $\omega_i$. We have $Q \subset P$ and $(P:Q) = 3$; let $\e: P \to \Z$/3$\Z$ be the homomorphism with kernel  $Q$ such that  $\e(\omega_1) = 1.$ [This choice of $\e$ means that we have ``oriented'' the Dynkin diagram, by choosing one of its extremities.]
  
  The scalar product on $P$ will be denoted by  $x\!\cdot\!y$; it takes values in ${1 \over 3}\Z$; for instance $\omega_1\!\cdot\!\omega_1={4 \over 3}$.
  We have $x\!\cdot\!y \in \Z$ if $x \in P, y\in Q$; if $\alpha$ is a root, we have $\alpha\!\cdot\!\alpha = 2$; we have $\alpha_i\!\cdot\!\omega_j = \delta^j_i$.\\

  Let now $L$ be the sublattice of \ $\Z \ \oplus $ $P$ made up of the pairs $(n,p)$ such that $n \equiv e(p) \pmod 3$. We define a scalar product
  $q_L$ on $L$ by the formula: $$ q_L(n,p;n',p') = nn'/3 - p\!\cdot\!p'.$$ 
  \hspace{3mm} Its values lie in $\Z$, and it is ``$\Z$-unimodular'', i.e., it gives an isomorphism of $L$ onto its $\Z$-dual. Its signature over $\R$ is (1,6); its discriminant is $1$.
  
   {\small[As we shall recall in \S7, the lattice $L$ is isomorphic to the N\'eron-Severi group of a smooth cubic surface, and the scalar product $q_L$ corresponds to the intersection form.]\par}
  
  \smallskip
  Note that  $Q$  embeds in $L$ by  $x \mapsto (0,x)$; this embedding transforms $q_L$ into the opposite of the scalar product of  $Q$.
  In particular, a root $\alpha$, viewed as an element of $L$, is such that  $q_L(\alpha, \alpha) = -2.$
  
  The intersection of $L$ with $\Z$ is generated by the element  $$h = (3,0).$$ We have  $q_L(h,h) = 3$ and $q_L(h,x) = 0$ if $x\in Q$.
  
  \n \ \ The element $\omega'_1 = (1,\omega_1)$ of $L$ is such that $q_L(\omega'_1,\omega'_1) = -1$ and $q_L(h,\omega'_1)~=1$; hence, if $\gamma = h-\omega'_1$, we have $q_L(\gamma,\gamma) = 0$.
  
  \smallskip
  
  {\small [From the cubic surface point of view of $\S7$, $h$ corresponds
  to a plane section, and $\omega'_1$ corresponds to a line. The formula $q_L(h,h) = 3$ means that the surface has degree 3 and the formula $q_L(h,\omega'_1)=1$ means that a line and a plane intersect in one point.]}
  
  \smallskip
  Let $G$ be the Weyl group of $R$, i.e., the subgroup of $\Aut(Q \otimes \R)$ generated by the reflections $s_\alpha$ associated with the roots $\alpha \in R$. The group $G$ acts on $P$ and $Q$. We extend its action to $\Z \ \oplus$ $P$, and hence to $L$, by making it act trivially on the factor $\Z$; the scalar product $q_L$
  is invariant by $G$.
  
  \medskip
  
   Proofs of the following theorem can be found in the standard texts on cubic surfaces (cf. [Ma 74], chap.IV, [De 80], [Do 12], chap.9).
   
   \smallskip
   
   \n {\bf Theorem 1.}  
   
   (a) {\it An element $\alpha$ of $L$ is a  root if and only if $q_L(h,\alpha)=0$ and
   $q_L(\alpha,\alpha)=-2.$}
   
   (b) {\it Let $Y$ be the set of $y \in L$ such that $q_L(h,y)=1$ and $q_L(y,y) = -1$. This set has $27$ elements, namely the pairs $(1,\omega)$ where $\omega$ belongs to the $G$-orbit of~$\omega_1$.}
   
   (c) {\it If $y,y' $ are two distinct elements of $Y$, then $q_L(y,y') = 0$ or $1$.}
  
  (d) {\it Let $\Omega$ be the graph with set of vertices $Y$, two vertices $y,y'$ being adjacent if and only if $q_L(y,y') = 1$. The
  natural injection $G \to \Aut(\Omega)$ is bijective.}
  
  \bigskip

  \n {\it Remark.}

  Since $\omega_1$ is orthogonal to the $\alpha_i$ for $i \geqslant 2$, it is fixed by the group  $H$ generated by $(s_{\alpha_2},...,s_{\alpha_6})$, which is a Weyl group of type $\sf{D}_5$, and has index 27 in $G$ since $|G| = 2^7\!\cdot\!3^4\!\cdot\! 5$ and $|H| = 2^7\!\cdot\!3\!\cdot\!5$. Hence $Y  \simeq G/H \simeq \Weyl(\sf{E}_6)/\Weyl(\sf{D}_5)$.

  \bigskip

{\bf \S2. The $27$-vertices graph $\Omega$ of Theorem $1$}.

\medskip

2.1. {\it Properties of $\Omega$.}

\smallskip
Every vertex belongs to 10 edges; every edge belongs to a unique triangle.
Hence the number of vertices, edges, triangles, tetrahedra is: 27, 135, 45, 0.

If  $x$ and $x'$ are two distinct and non adjacent vertices, the element $x-x'$
of $L$
is a root: this follows from Theorem 1 (a) since $q_L(x-x',x-x')= -1-1=-2$ and $q_L(h,x-x')=0$. Every root is obtained in such a way by $six$ disjoint couples
$(x_i,x'_i)$, $ i= 1,...,6$. Every $x_i$ is adjacent to every $x'_j$ with $j\ne i$.
The reflection defined by the root $\alpha = x_1-x'_1=...=x_6-x'_6$ permutes
$x_i$ and $x'_i$ for every $i$, and fixes the other 15 vertices of $\Omega$.
The twelve vertices $x_1,...,x_6,x'_1,..., x'_6$ make up a {\it double-six}, see $\S$2.2 below.

Let $g$ be an element of $G$ of order 2. Assume it is not a reflection. The number
of the vertices fixed by $g$ is then either 7 or 3. The first case occurs if and only if
$g$ is a product of two reflections.

\medskip  

2.2. {\it An explicit construction of $\Omega$.}

  Let us explain how one can describe the graph  $\Omega$ in combinatorial terms, following Schl\"afli [Sch 58]. 
  
  Let $X = \{1,2,...,6\}$ and let
 $X'$ be a copy of  $X$; if $x \in X$, the corresponding
  point of $X'$ is denoted by $x'$; let $S$ be the set of all subsets of $X$ with 2 elements. The graph 
  $\Omega$ of Theorem 1 is isomorphic to the graph $\Omega_X$ whose set of vertices $Y$ is the disjoint union $X \sqcup X' \sqcup S$, two vertices being adjacent
  in the following cases (and only in those):
  
  \smallskip
  
    $x \in X$ \ adjacent to $y' \in X' \  \ \Longleftrightarrow \ \ x \neq y$,
    
    $x \in X$ \ adjacent to  $s \in S \ \ \ \ \Longleftrightarrow \ \ \ x \in s$,
    
    $x' \in X'$ adjacent to  $s \in S \ \ \ \Longleftrightarrow\ \ \ x \in s$,
    
    $s_1 \in S$ \ adjacent to $s_2 \in S \ \ \Longleftrightarrow s_1 \cap s_2 = \varnothing$.
    
    \medskip 
    
    {\small Note that every element $\{x,y\}$ of $S$ defines two edges: $x$\! --\! $y'$
    and $y$\! -- \!$x'$. There are no edges connecting vertices that are both in $X$, or both in $X'$: the full subgraph of $\Omega$ with set of vertices the double-six $X \sqcup X'$ is a
    bipartite graph.}
   
   \smallskip
   The group $\sf{S}_6$ of permutations of $X$ acts on $\Omega_X$ by {\it transport de structure}; an element of that group is a reflection in $G$ if and only if it is a
   transposition.
  Let $\varepsilon$ be the automorphism of order $2$ of  $\Omega_X$ which fixes the points of $S$ and exchanges $x\in X$ with $x' \in X'$;  that automorphism is also a reflection and it commutes with the action of $\sf{S}_6$. We thus obtain an embedding of the group $\{1,\varepsilon\} \times \sf{S}_6$ into $\Aut(\Omega_X)$.
   From the Weyl group point of view, this corresponds to the embedding of $\Weyl(\sf{A}_1 \times \sf{A}_5) \simeq \{1,\varepsilon\} \times \sf{S}_6 $ into  $\Weyl(\sf{E}_6)$
  defined by the inclusion $\sf{A}_1 \times \sf{A}_5  \to  \sf{E}_6$.
  
\bigskip

{\bf \S3. Involutions and cubes in $ G = \Weyl(\sf{E}_6)$}

\medskip
3.1. {\it Finite Coxeter groups.}

\smallskip

We recall first a few definitions and results that apply to every finite Coxeter group $G$, and are easy to prove.

Let $g \in G$ be  an {\it involution}, i.e. an element $g$ such that $g^2=1$.
The multiplicity of $-1$, as an eigenvalue of $g$ in the Coxeter representation, is called the
 {\it degree} of $g$. An involution of degree 1 is a reflection. 

A {\it cube} (cf. [Se 18]) of $G$ is an abelian subgroup of $G$ generated by
reflections. It is elementary abelian of type (2,...,2). If $C$ is a cube, the
set $S_C$ of the reflections contained in $C$ is a basis of $C$ (viewed as an $\F$$_2$-vector space), i.e. every element $g$ of $C$ can be written uniquely as
$g = \prod_{s\in I} s$, where $I$ is a subset of $S_C$, and we then have $ \deg(g)=|I|$. If the rank
of $C$ is $n$, i.e. if $|S_C|=n$, then $C$ contains a unique element of degree
$n$, namely the product of all the $s\in S_C$; we call that element
the {\it extremity} of $C$.

Every involution is the extremity of at least one cube.

From the root system point of view, a reflection corresponds to a pair of opposite
roots, and a cube corresponds to a family of mutually orthogonal roots.

\medskip

The following properties are valid for every Weyl group $G$, and they can easily be checked
in the special case where $G$ is of type $\sf{E}_6$, see $\S$3.2 below.
 
\smallskip

(3.1) {\it The involutions of $G$ of maximal degree are $G$-conjugate.}

\smallskip

(3.2) {\it If $G$ is simply-laced~\footnote{Recall that a Coxeter group is called {\it simply laced} if it is a product of Weyl groups of types $\sf{A}, \sf{D}, \sf{E}$;
this is equivalent to asking that the product of two reflections has order 1, 2 or~3.}, 
any two cubes of $G$ with the same extremity are $G$-conjugate.}

\smallskip

(3.3) {\it If $C$ is a maximal cube of $G$, the centralizer of $C$ in $G$ is  equal to $C$.}
 
 \medskip
 
 3.2. {\it The involutions of $ \Weyl(\sf{E}_6)$.}

\medskip

\n From now on, we assume again that $G$  is of type $\sf{E}_6$, as in $\S$1. Then:

\smallskip

(3.4) {\it The involutions have degree $0, 1, 2, 3$ or $4$.}

\smallskip

(3.5) {\it Any two involutions of the same degree are $G$-conjugate.} 

\smallskip
(3.6) {\it The number of involutions of degree $0, 1, 2, 3, 4$ is $1, 36, 270, 540, 45$ respectively.}

\smallskip

These three facts can be read off from the character table of $G$ (see [ATLAS], pp.26-27, where our $G$ is denoted by $G.2$). Another possible proof of (3.4) is to use the fact that every involution 
is contained in the centralizer  of a reflection; that centralizer
is a Weyl group of type $\sf A_1 \times \sf A_5$, and we only have to check that
the involutions of $\Weyl(\sf A_5) \simeq \sf{S}_6$ have degree 0, 1, 2 or 3. The same argument proves (3.2) and (3.3), and also (3.9)  below.

\medskip

3.3. {\it The maximal cubes.}

\smallskip 

\n By (3.3) and (3.4), we have

(3.7) {\it The maximal cubes of $G$ have order $2^4$, and they are $G$-conjugate.}

\smallskip

\n This implies:

(3.8) {\it A maximal cube of $G$ contains $1, 4, 6, 4, 1$ involutions of degrees $0, 1, 2, 3, 4$ respectively.}

\smallskip

(3.9) {\it The number of maximal cubes is $135$.}

\smallskip
\n {\it Examples.}

a) The roots  $(\alpha_2,\alpha_3,\alpha_5, \alpha_2+\alpha_3+\alpha_5+2\alpha_4)$ are mutually orthogonal. Hence they define
a maximal cube. That cube is contained in the $\sf{D}_4$-subroot system of $\sf{E}_6$ with basis $ (\alpha_2,\alpha_3,\alpha_4,\alpha_5)$. Another example
  is the cube defined by $(\alpha_1,\alpha_4, \alpha_6, \alpha_0)$.
    
 b) In terms of the combinatorial description of \S2, we may choose as maximal cube the group generated by the following four reflections: the automorphism $\varepsilon$ and the three transpositions (12),(34),(56) in $\sf{S}_6$.

  \medskip
  
\n  {\it The normalizer of a maximal cube}.
  
  \smallskip
  
    Let $C$ be a maximal cube and let $N$ be its normalizer in $G$. The group
    $N/C$ acts on $C$, and hence permutes the four reflections of $C$.
    We thus get a homomorphism $N/C \to \sf{S}_4$. 
  
   \smallskip

  (3.10) {\it The homomorphism $N/C \to \sf{S}_4$ is bijective.}

   \smallskip

The injectivity follows from (3.3). The surjectivity can be proved, either by
writing explicitly enough elements of $N$, or by counting: since the number of cubes
is 135, and the order of $G$ is 51840, the order of $N$ is $51840/135=384$ and the order of $N/C$ is $384/16=24$.

\medskip
  
  3.4. {\it Action of a maximal cube on the vertices of the graph $\Omega$.} 
  
  \smallskip
  
  Let $C$ be a maximal cube of $G$. We shall later need  the following information on the action of $C$ on the set $Y$ of Theorem 1 (the ``twenty-seven lines").

\medskip

\ \n {\bf Lemma 1.} {\it The action of $C$ on $Y$ has three fixed points, and six orbits of order $4;$ these orbits are isomorphic to $C/C_i,
\ i = 1,...,6$, where the $C_i$ are the six cubes of order $4$ contained in $C$.}

\medskip

\n {\it Proof}. We use the description of $Y$ of $\S$2.2, as $Y= X \sqcup X' \sqcup S$,
and we use for $C$ the cube defined in Example b) of $\S$3.3, namely the one
generated by the reflections $\varepsilon, (12), (34), (56)$.

The action of $C$ fixes the points $\{1,2\}, \{3,4\}, \{5,6\}$ of $S$. The four points $\{1,3\}, \{1,4\}, \{2,3\}, \{2,4\}$ of $S$ make up an orbit
isomorphic to $C/C_1$ where $C_1$ is the subgroup of $C$ generated by $\varepsilon$ and the transposition (56); there are two other similar orbits in $S$.
In $X \sqcup X'$, the points $1, 1', 2, 2'$ make up an orbit isomorphic to $C/C_4$, where $C_4$ is generated by the transpositions (34) and (56); there are two other similar orbits.

\medskip

{\bf \S4. The $G$-quadratic forms $q_6, q_7, q_{27}$ and the $C$-quadratic form} $q_4$.

\medskip

4.1.   Let $G$ be a group. A {\it $G$-bilinear form} over a commutative ring  $A$ is a bilinear form $q$ over an
  $A$-module $E$, together with an $A$-linear action of $G$ on $E$ that fixes  $q$. 
  
  In the above case, where
  $G = \Weyl(\sf{E}_6)$, we have three such examples, with $A = \Z$, which are of rank 6, 7 and 27.
  
  \smallskip
  (i) The bilinear form  $\alpha\!\cdot\!\beta$ on  $Q$, cf. $\S1$.
  
  \smallskip
  (ii) The bilinear form  $q_L$ on $L$.
  
  \smallskip
  (iii) The bilinear form $q_Y$ on $\Z$$^Y$, given by $q_Y(e_y,e_{y'}) = \delta_y^{y'}$, where  $(e_y)_{y \in Y}$ is the natural basis of $\Z$$^Y$.
  
  \smallskip
  In each case, the action of $G$ is the obvious one, and the form is symmetric.
  
 \medskip
 
 4.2. In what follows, $k$  is a field of characteristic $\neq 2$. This will allow us to identify quadratic forms and symmetric bilinear forms, hence to write  $q(x)$
 instead of $q(x,x)$.
 
 The base change $\Z \to$ $ k$ transforms the three $G$-bilinear forms above 
 into $G$-quadratic forms over $k$. 
We denote them by $q_6, q_7, q_{27}$, but, for convenience, we divide by $2$ the first one. In other words, we put 
 
 \smallskip
(4.1)  $q_6(\alpha, \beta) = {1 \over 2}\alpha.\beta$, 

\smallskip
\n  so that  $q_6(\alpha,\alpha) = 1$ for every $\alpha \in R$.
 
 \smallskip
  The forms $q_7$ and $q_{27}$ are non-degenerate. The same is true of $q_6$
   provided $\ch(k)\ne\!3$.

\smallskip 
  With standard notation, we have an isomorphism of $G$-quadratic forms:
 
 \smallskip
 (4.2) $  q_7 = \langle 3 \rangle + \langle -2 \rangle q_6$ \ {\it if} $\ch(k) \ne 3$,
 
 \smallskip
\n with $G$ acting trivially on the rank $1$ form $\langle 3 \rangle: $ this follows from
the definition of $q_L$ in $\S$1.

 \smallskip
    There is no such simple formula for the $G$-form $q_{27}$, since the underlying linear representation of $G$ is not a linear combination of exterior powers of the standard degree 6 representation.
  
  \medskip
  
  4.3.  {\it The $C$-quadratic form $q_4$.}
  
  \smallskip 
  Let $C$ be the maximal cube introduced in \S3 from the combinatorial point of view. It contains four reflections $s_1,...,s_4$, namely $\varepsilon$ and the transpositions (12), (34) and (56). Let 
  $\beta_1,...,\beta_4$ be the roots (well-defined up to signs) corresponding to $s_1,...,s_4$. Let $V$ be the $k$-vector subspace of $Q  \otimes k$ generated by the $\beta_i$, and let $q_4$ denote the restriction of $q_6$ to $V$; we have $q_4(\beta_i,\beta_j) = \delta_i^j$. The group $C$ acts in a natural way on $V$; hence $q_4$ {\it is a $C$-quadratic form of rank} 4. The space $V$ splits under the action of $C$ into four 1-dimensional subspaces, orthogonal to each other. This gives a splitting of $q_4$ as 
  
  \smallskip
  
  (4.3)  $q_4 = r_1 + \cdots + r_4$,
  
  \smallskip
  
  \n where $r_i$ is the $C$-quadratic form of rank $1$ generated by $\beta_i$, on which  $s_i$ acts by $-1$ and the other $s_j$ act trivially.
  
  \medskip
  
  4.4. {\it Relations between the $G$-quadratic forms $q_6, q_7, q_{27}$ and the $C$-quadratic form} $q_4$.
  
    Any $G$-quadratic form $q$ defines, by restriction of the action of the group, a $C$-quadratic form, which we denote by $q|C$. This applies
    in particular to $q_6, q_7, q_{27}$. The $C$-quadratic forms so obtained can all be expressed in terms of $q_4$. Namely
    
    \medskip
    
    \n {\bf Theorem 2.} {\it We have the following isomorphisms of $C$-quadratic forms} 
    
    \smallskip
    (4.4)\ $q_6|C \ = \ q_4 + \langle 2,6 \rangle$ \quad {\it if} \ $\ch(k)\ne 3$,
    
    \smallskip
  (4.5)\ $q_7|C \ = \ \langle-2\rangle q_4 + \langle -1,-1,1 \rangle,$
  
  \smallskip
  (4.6)\ $q_{27}|C = \lambda^2q_4 + 3q_4 + 9.$
  
  \medskip
  
 \n [Here, $\langle 2,6 \rangle$ denotes the quadratic form  $\langle 2 \rangle +  \langle 6 \rangle$ with trivial action of $C$. Similarly, $3q_4$ means
   $\langle 1,1,1 \rangle \otimes q_4= q_4 + q_4 + q_4$ and 9 means the direct sum of nine copies of $\langle 1 \rangle$ with trivial action of $C$. As for $\lambda^2q_4$, it is the second exterior power
  of $q_4$, with the natural action of $C$, cf. [Se 03], $\S$27.1.]
  
  \medskip
  
  \n {\it Proof of} (4.4). This is a simple computation in the root lattice $Q$. With Bourbaki's notation ([Bo 68], \S4.12), we may choose for the $\beta_i$ of $\S$4.3 the roots $\varepsilon_1 + \varepsilon_2$, $\varepsilon_1 - \varepsilon_2$, $\varepsilon_3 + \varepsilon_4$, $\varepsilon_3 - \varepsilon_4$. They are orthogonal
  to $\gamma = \varepsilon_5$ and $\delta = \varepsilon_8 - \varepsilon_6 - \varepsilon_7$. Moreover, $\gamma$ and $\delta$ are orthogonal to each other,
  ${1\over 2}\gamma.\gamma = {1\over 2}$, and ${1\over 2}\delta.\delta = {3\over 2}$. Hence the orthogonal complement $V'$ of $V$ in $Q \otimes k$ is quadratically isomorphic to $\langle {1\over 2},{3\over 2} \rangle = \langle 2,6 \rangle$, and the action of $C$ on $V'$ is trivial.

  \medskip
  
  \n {\it Proof of} (4.5). Let $V''$ be the orthogonal complement of $V$ in $L \otimes k$. The restriction $q_3$ of $q_L$ to $V''$ is a rank 3 quadratic form. We have $q_7|C = q_4|C +  q_3$. The action of $C$ on $q_3$ is trivial: indeed, the reflections which generate $C$ fix every element orthogonal to the corresponding roots, hence fix $V''$. Since the discriminants of $q_7$ and $q_4$ are equal to 1,
  the same is true of $q_3$. By Lemma 1, $C$ fixes three points
  of $Y$. Let $y$ be one of them, and let $x = h-y$. Both $x$ and $y$
  belong to $V''$. Since $q_7(h,h)=3, q_7(h,y)=1$  and  $q_7(y,y)=-1$, we have $q_L(x,x) = 0$; moreover, $x$ is nonzero since $q_7(x,h)=2$. Hence $q_3$ is isotropic; it contains the 2-dimensional quadratic form $\langle-1,1\rangle$. Since its discriminant is $1$, we have  $q_3=\langle-1,-1,1\rangle$, as claimed.
  
  \medskip
    
  \n {\it Proof of} (4.6). Lemma 1 gives a decomposition of $q_{27}|C$ as the orthogonal sum of $ \langle 1,1,1 \rangle$ with trivial action, and six 
  $C$-quadratic forms $q'_i$ associated with the permutation sets $C/C_i$, where the $C_i$ are the six cubes of order 4 contained in $C$. Consider for instance the case where $C_i$ is generated by $s_1,s_2$, as in \S4.3. In that case, one finds that $q'_i = 1 + r_3 + r_4 + r_3r_4$, where the $r_i$ are the $C$-quadratic forms of rank 1 occurring in (4.3). The other $q'_i$ correspond similarly to the pairs $\{1,3\}, \{1,4\}, \{2,3\}, \{2,4\}, \{3,4\}$. Adding up gives 
  
  \medskip
  
  $ q_{27}|C = 3 + 6 + 3 \ \Sigma_{n=1}^4 \ r_n + \Sigma_{1 \leqslant  m < n \leqslant 4}\ r_mr_n.$
  
  \medskip
  
 \n  By (4.3), we have $ \Sigma_{n=1}^{n=4} \ r_n = q_4$ and $\Sigma_{1 \leqslant  m < n \leqslant 4}\  r_mr_n = \lambda^2q_4$. We thus obtain (4.6).

  \medskip
  
  \n {\it Remark}. One may also prove (4.4) by combining (4.5) and (4.2). 
  
\bigskip
{\bf \S5. Twists.}

\medskip
5.1. We keep the assumption that the characteristic of $k$ is $\neq 2$.
  Let $k_s$ be a separable closure of $k$. Let $\Gamma_k  = \Gal(k_s/k )$ be the ``absolute Galois group'' of $k$.
  If $\varphi: \Gamma_k \to G$ is a continuous homomorphism, we may view $\varphi$ as a 1-cocycle of $\Gamma_k$ with values in the orthogonal group $\Or_q(k_s)$, and use it to {\it twist}   $q$,  cf. [Se 64], chap.III, $\S$1.2. We thus get a quadratic form $q_\varphi$ over $k$. 
  
  \medskip
  \n {\small [Let us recall an explicit construction of $q_\varphi$.
  
 Let $V$ be the $k$-vector space on which $q$ is defined and $G$ acts.  On
 $V_s = k_s \otimes V$ there is  a
 natural $k_s$-linear action of $G$ (denoted by  $gx$ if $g\in G$ and $x\in V_s$)  and also a natural semi-linear action of $\Gamma_k$ (denoted by $\gamma x$ if $\gamma \in \Gamma_k$ and $x \in V_s$).
 These two actions commute. Hence they can be combined to define another semi-linear action  of $\Gamma_k$ on $V_s$, namely $\gamma\bullet x = \gamma(\varphi(\gamma)x).$

  Let
 $V_\varphi$ be the set of $x\in V_s$ such that $\gamma\bullet x=x$ for every
 $\gamma \in \Gamma_k$. It is a $k$-vector subspace of $V_s$, and the natural map 
 $k_s \otimes V_\varphi \to V_s$ is an isomorphism,  i.e. $V_s$ is a ``$k$-form''
 of $V$, cf. [Bo 81], $\S$10, prop.7. The form $q$ extends to a $k_s$-quadratic form $q_s$ of $V_s$, which is $G$-invariant and $\Gamma_k$-equivariant: $q_s(\gamma x) = \gamma q_s(x)$ for every
 $\gamma \in \Gamma_k$ and every $x \in V_s$. 
 
   If $x \in V_\varphi$, then, for every $\gamma \in \Gamma_k$, we have $x = \gamma(\varphi(\gamma)x)$, hence:
   $$q_s(x) = \gamma(q_s(\varphi(\gamma)x)) = \gamma(q_s(x)),$$
      which shows that $q_s(x)$ belongs to $k$. Hence the restriction of $q_s$ to $V_\varphi$
   is a $k$-valued quadratic form on $V_\varphi$: it is the quadratic form $q_\varphi$
   we wanted to describe.
 
  \smallskip
  \n {\it Remark}. Let $x\in V$ be fixed under the action of $\varphi(\Gamma_k)$. Then $x$ belongs to $V_\varphi$ and  $q_\varphi(x) = q(x)$. In particular:
   
   \smallskip
  \n  (Iso) {\it If $V$ contains a non-zero isotropic vector fixed under $\varphi(\Gamma_k)$, then $q_\varphi$
   is isotropic.}}]
   
   \medskip
  
5.2.  This construction applies in particular to the $G$-forms  $q_6, q_7, q_{27}$ above; we thus obtain the quadratic forms $q_{6,\varphi}, q_{7,\varphi}$ and $q_{27,\varphi}$. The forms $q_{7,\varphi}$ and $q_{27,\varphi}$ are non-degenerate. The form $q_{6,\varphi}$ is non-degenerate if $\ch(k)\ne 3$.

\smallskip
  
  \n{\bf Theorem 3.} {\it The form $q_{7,\varphi}$ is isotropic}.  
   
  \smallskip
  
  \n {\it Proof.} Assume first that $\Gamma_k$ fixes a point $y$ of $Y$. Then it also fixes $x= h-y$, and we have $q_7(x)=q_7(h)-2q_7(h,y)+q_7(y)= 3-2-1=0.$ By (Iso) above this shows that $q_{7,\varphi}$ is isotropic.
  
   In the general case, since the set $Y$
  has odd order, the same is true of one of the orbits of $\Gamma_k$, i.e., there exists $y\in Y$ which is fixed by an open subgroup of $\Gamma_k$ of odd index.
  That subgroup is of the form $\Gamma_{k'}$
with $k'/k$ of odd degree. Hence $q_{7,\varphi}$ becomes isotropic after the base change $k \to k'$. By a well-known theorem of Springer ([Sp 52]), this implies that $q_{7,\varphi}$ is isotropic over $k$.

\medskip

\n {\bf Corollary 1.} {\it There exists a unique quadratic form $q_{5,\varphi}$ over $k$, of rank $5$, such that} 

\smallskip

$(5.1) \ q_{7,\varphi} = \langle-2\rangle q_{5,\varphi} + \langle1,-1\rangle.$

\smallskip
\n {\small (The factor $\langle-2\rangle$ is included in order to have formula (5.2) below.)}

\medskip

\n {\bf Corollary 2.} {\it If $\ch(k) \ne 3, $ we have}:

\smallskip
$(5.2) \ q_{6,\varphi} = q_{5,\varphi} + \langle6\rangle.$

\smallskip
$(5.3) \ q_{7,\varphi} = \langle-2\rangle q_{6,\varphi} + \langle3\rangle.$

\medskip

\n Formula (5.3) follows from (4.2). By combining it with (5.1), we get 

\smallskip
(5.4) \  $\langle-2\rangle q_{6,\varphi} + \langle3\rangle = \langle-2\rangle q_{5,\varphi}  +
  \langle3,-3\rangle$, \ {\it since} $\langle1,-1\rangle = \langle3,-3\rangle$,
  
  \smallskip
  \n hence $\langle-2\rangle q_{6,\varphi} = \langle-2\rangle q_{5,\varphi} + \langle-3\rangle$, which is equivalent to (5.2).

\medskip 
\n {\it Remark}. When $\ch(k)\ne 3$, the quadratic form $q_{5,\varphi}$ is {\bf not}
the twist of a $G$-quadratic form of rank 5, if only because $G$ has no non-trivial
linear representation of degree $5$. 

\bigskip

5.3. {\it  The   quadratic form $q_{27,\varphi}$.}

\smallskip

This form may be viewed as a {\it trace form}. Indeed, the group $\Gamma_k$ acts on $Y$ via 
    $\varphi$, and this defines an \'etale algebra of rank 27 over $k$ whose trace form is $q_{27,\varphi} $, cf. e.g., [BS 94], \S1.
    
    It can be expressed
in terms of $q_{5,\varphi}$, namely:
    
    \bigskip
 
 \n {\bf Theorem 4.} {\it We have} 
 
 \smallskip 
 (5.5) \  $q_{27,\varphi} =    \lambda^2q_{5,\varphi} \ + \ \langle1,2\rangle q_{5,\varphi} \ + 6 +\langle2\rangle .$
 
 \medskip
 
   The proof will be given in \S6, using the method of [Se 03] (extended in [Se 18] to all Weyl groups): checking first the case where $\varphi: \Gamma_k \to G$ takes values in a maximal cube, and then showing that this special case implies the general one.

\medskip

\n  {\it Remark.}

Here  also, formula (5.5) can be rewritten in terms of $q_{6,\varphi}$  and $q_{7,\varphi}$. One finds 

\smallskip

(5.6)  $q_{27,\varphi} =    \lambda^2q_{6,\varphi} \ + \ \langle 3 \rangle q_{6,\varphi} \ + \ 6,$ \ \ if $\ch(k) \ne 3$.

\n and 

(5.7) $q_{27,\varphi} =    \lambda^2q_{7,\varphi} \ + \ (\langle -1 \rangle - \langle2\rangle)q_{7,\varphi} \ + \ 7 - \langle-2\rangle.$

    \bigskip
{\bf \S6. Proof of Theorem 4.}

\medskip
\n 6.1. {\it Proof of Theorem 4 when $\varphi$ maps $\Gamma_k$ into the cube }$C$.

\smallskip

In that case, we may use $\varphi$ to twist the $C$-quadratic form $q_4$ of $\S$4.3, and we obtain a quadratic form $q_{4,\varphi}$. Formulas (4.5) and (4.6) imply 

\smallskip
(6.1) \ $q_{7,\varphi} = \langle-2\rangle q_{4,\varphi} + \langle-1,-1,1\rangle$.

\n and
\smallskip

(6.2) \ $q_{27,\varphi} = \lambda^2q_{4,\varphi} +3q_{4,\varphi} + 9$.

\smallskip
\n Since $3 = \langle1,1,1\rangle = \langle1,2,2\rangle$, we may rewrite (6.2) as 

\smallskip
(6.3) \ $q_{27,\varphi} = \lambda^2q_{4,\varphi} +\langle1,2,2\rangle q_{4,\varphi} + 9$,

\medskip

\n By (5.1), we have  

\smallskip
(6.4) $q_{7,\varphi} = \langle-2\rangle q_{5,\varphi} + \langle1,-1\rangle.$ 

\smallskip

\n Comparing (6.1) and (6.4) gives

\smallskip

(6.5) \ $q_{5,\varphi} = q_{4,\varphi} + \langle2\rangle,$

\smallskip
\n hence

\smallskip

(6.6) \ $\lambda^2q_{5,\varphi} = \lambda^2q_{4,\varphi}+ \langle2\rangle q_{4,\varphi}.$

\smallskip

\n Using (6.4) and (6.5) we may rewrite (6.3) as 

\smallskip

(6.7) \ \ \ $q_{27,\varphi} = \lambda^2q_{5,\varphi} - \langle2\rangle q_{4,\varphi} + \langle1,2,2\rangle q_{4,\varphi} +9 $

\hspace{19mm} $= \lambda^2q_{5,\varphi} +\langle1,2\rangle q_{4,\varphi} + 9$

\hspace{19mm} $= \lambda^2q_{5,\varphi} +\langle1,2\rangle(q_{5,\varphi}-\langle2\rangle)+9$

\hspace{19mm} $= \lambda^2q_{5,\varphi} +\langle1,2\rangle q_{5,\varphi} + 9 - \langle1,2\rangle$

\hspace{19mm} $= \lambda^2q_{5,\varphi} +\langle1,2\rangle q_{5,\varphi} + 6 + \langle2\rangle$,

\smallskip

\n which is (5.5).

\medskip

    \n 6.2. {\it The Witt-Grothendieck invariants defined by $ q_4, q_5, q_6, q_7, q_{27}$.}
    
    \smallskip
     (In what follows, we use freely the definitions and the elementary properties of the ``invariants" of a finite group given in the first sections of [Se 03].)
     
     \smallskip
    
    The cohomology set $H^1(k,G)$ of all $G$-torsors over $k$ may be identified with 
    the conjugation classes of continuous homomorphisms $\varphi: \Gamma_k \to G$, cf. e.g. [BS 94], \S1. Since
    the Galois twists defined by conjugate homomorphisms are isomorphic, we may interpret the maps $\varphi \mapsto q_{5,\varphi},
 q_{7,\varphi}, q_{27,\varphi}$ as maps of $H^1(k,G)$ into the Witt-Grothendieck ring $WGr(k)$ of $k$; let us denote them by $a_{5,k}, a_{7,k}, a_{27,k}$.
 This construction applies to all the field extensions $K$ of $k$, and we thus obtain three Witt-Grothendieck invariants  $a_5,a_7, a_{27}$ of $G$,
 i.e. three elements of the group $\Inv_k(G,WGr)$, cf. [Se 03], $\S$27.1. 
 We also get an invariant $a_4$ by $a_4 = a_5 - \langle2\rangle$; note that
 this invariant is not in general {\it effective}: its values cannot always be
 represented by quadratic forms of  degree 4 (but they can be if  $\varphi$ takes values in $C$). Moreover, when $\ch(k)\ne 3$, we have the invariant $a_6 = a_7 - \langle6\rangle$.
 
 In what follows, we  mainly use $a_4$.
 
 \medskip
 
  For every subgroup $H$ of $G$, there is a
 natural restriction map $$\Inv_k(G,WGr) \to \Inv_k(H,WGr),$$ cf. [Se 03], \S13.
 
 \medskip
 
  \n 6.3. {\it The splitting principle.}

 \medskip
 
   A basic fact about such invariants is:
   
  \medskip
   
   \n {\bf Theorem 5.} {\it Let $\cG$ be a Weyl group, and let $a \in \Inv_k(\cG,WGr)$. Assume that
   the restriction of $a$ to every cube of $\cG$ is $0$. Then $a=0$}.
   
   \medskip

This is proved in [Se 03], \S29, when $\cG$ is a symmetric group (for Witt invariants - the case of the Witt-Grothendieck invariants follows). The
proof for an arbitrary Weyl group is similar;  it has not been published yet,
but we hope it will soon be.

\smallskip

  In the Appendix, we shall need the analogue of Theorem 5 for the case of the
  group $\Inv_k(\cG)$ of the {\it cohomological invariants} mod 2 of $\cG$ (i.e.
  the group $\Inv_k(\cG,\Z/2\Z)$ of [Se 03], $\S$4.1):
  
  \medskip
  
\n {\bf Theorem 6.} {\it Let $a$ be an element of $ \Inv_k(\cG)$. If the restriction
of $a$ to every cube of $\cG$ is $0$, then $a=0$.}

\medskip
A proof can be found in [GH 19], under the assumption that $\ch(k)$ does not
divide $|\cG|$; that assumption (which, for $\Weyl(\sf{E}_6)$, would eliminate characteristics 3 and 5) is not necessary, as one sees by the method sketched in [Se 18].

   \medskip
   
 \n 6.4. {\it End of the proof of Theorem 4.}

   \smallskip
   We apply Theorem 5 with $\cG = G$, and with $a \in \Inv_k(G,WGr)$ defined by 
   
   \smallskip
   
   (6.8) $a = a_{27} -  \lambda^2a_5  - \langle 1,2 \rangle a_5 - 6 - \langle2\rangle.$
   
   \smallskip

By \S6.1, the restriction of $a$ to $C$ is $0$. Since every cube $C'$ of $G$ is conjugate to a subgroup of $C$, the restriction of $a$ to $C'$ is 0.
By Theorem 4, this implies $a = 0$; hence (5.5).

\smallskip

\n {\it Remarks.}

  1.
The same method can be used to describe the structure of $\Inv_k(G,WGr)$. The result is simpler to state for the Witt invariant ring 
$\Inv_k(G,W)$: this ring is a free $W(k)$-module with basis the five elements $\lambda^i a_4,  i = 0,...,4$ (or, equivalently, the $\lambda^ia_5$ or the $\lambda^ia_7, i=0,...,4).$

\smallskip
    
    2. Let $T$ be a finite $G$-set, and let $\varphi: \Gamma_k \to G$, be as above;
    we thus have an action of $\Gamma_k$ on $T$, hence 
 an \'etale $k$-algebra $E_{T,\varphi}$.
    (For the dictionary ``finite $\Gamma_k$ -set $\Longleftrightarrow$ \'etale $k$-algebras'', see [Bo 68], chap.V, $\S$10.10 and [KMRT 98], $\S$18.A.)

     Let
    $q_{T,\varphi}$ be its trace form. It follows from Remark 1 above that we have
     $$ (6.9) \ \ q_{T,\varphi} = t_0 + t_1q_{4,\varphi} + t_2\lambda^2q_{4,\varphi}+t_3\lambda^3q_{4,\varphi} + t_4\lambda^4q_{4,\varphi},  \hspace{20mm}$$
     \n   where the $t_i$ are quadratic forms of type either $\langle1,...,1,1\rangle$ or
    $\langle1,...,1,2\rangle$. (Note that $\langle2,2\rangle = \langle1,1\rangle$, so that
    there is no point in having more $\langle2\rangle$'s.)
    All what is needed is the decomposition of $T$ with respect to the action
    of $C$; the orbits of order $2$ or $8$ are the ones which introduce a
    factor $\langle2\rangle$.
    
    \smallskip
    
      A typical example is when $T$ is the set of the 45 triangles of $\Omega$
      (they correspond to the tritangent planes in the cubic surface interpretation of $\S7$). One finds that $C$ has one fixed point, six orbits of order 2 and four orbits of order 8. By looking at the structure of these orbits, one gets:
      
      \smallskip
      $ (6.10) \ \  q_{T,\varphi} = 11 + \langle1,1,2\rangle q_{4,\varphi} + \langle1,1,2\rangle\lambda^2q_{4,\varphi}+\langle2\rangle\lambda^3q_{4,\varphi}.$
      
      \smallskip
 \n This formula becomes slightly simpler when expressed in terms of $q_{5,\varphi}$:
 
 \smallskip
 
 $ (6.11)    \ \  q_{T,\varphi} = 9 + \langle2\rangle + q_{5,\varphi} + \langle1,2\rangle\lambda^2q_{5,\varphi}+\langle2\rangle\lambda^3q_{5,\varphi}.$

    \bigskip
    
    {\bf \S7. The cubic surfaces and their $27$ lines.}

\medskip

\n7.1. In this subsection, we drop our assumptions on the ground field $k$, i.e., we allow char$(k) = 2$.

\smallskip
Let $V \subset {\bf P}_3$ be a smooth cubic surface over $k$. It is well known that, over a suitable field extension of $k$, it
contains 27 lines, cf. [Ha 77], chap.5.4, [Do 12], chap.9  and [Ma 74], chap.IV. These lines are rational  over $k_s$, cf. [Co 88], Theorem 1 and [KW 17], Corollary 53.
Let $L_V$ be the N\'eron-Severi group of $V$ over $k_s$, equipped with the symmetric $\Z$-bilinear form ``intersection product''.
It is a lattice of rank 7, and it contains the following elements 

\smallskip
  (a) The class $h_V$ of the hyperplane sections of $V$; we have $h_V\!\cdot \!h_V = 3$.
  
  \smallskip
  (b) The set $Y_V$ of the classes of the 27 lines; if $y \in Y_V$, we have $y\!\cdot \!y = -1$ and $h_V\!\cdot \!y$ = 1; if $y' \in Y_V$ is distinct from
  $y$, we have $y\!\cdot \!y' = 0$ if the corresponding lines are disjoint, and $y\!\cdot \!y' = 1$ if they meet.
      
      \smallskip
      
      It is well known that the triple $T_V = (L_V, h_V,$ intersection product) is isomorphic to the triple $T =(L, h, q_L)$
      of \S1; see the above  references.  More precisely, let $\Theta_V$ denote the set of isomorphisms $\theta: T_V \to T$.
      We have a left action of $G = \Weyl(\sf{E}_6) = \Aut(T)$ on $\Theta_V$, by $g\theta = g\circ \theta$; that action is free, and transitive. On the other hand,
    we have a right action of $\Gamma_k$ on $\Theta_V$, by $\theta \gamma=  \theta \circ \gamma$, for $\gamma \in \Gamma_k$. These two actions commute. We may view such a situation in the following equivalent ways (cf. [BS 94], \S1.3):
    
    \smallskip
    (7.1) The action of $\Gamma_k$ on $\Theta_V$ defines an \'etale algebra $E_V$, on which $G$ acts, and one hence gets a {\it $G$-Galois algebra
    over} $k$.
    
    (7.2) The $k$-finite \'etale scheme Spec$ \ E_V$ is a $G$-torsor over $k$, whose set of $k_s$-points is $\Theta_V$.
    
    (7.3) If we choose a point $\theta$ of $\Theta_V$, for every $\gamma \in  \Gamma_k$ there is a unique element $\varphi(\gamma)$ of $G$
    such that $\varphi(\gamma)\theta = \theta\gamma$, and the map $\gamma \mapsto \varphi(\gamma)$ is a continuous homomorphism  $\varphi: \Gamma_k \to G$. Changing the choice of $\theta$ replaces $\varphi$ by a conjugate.
    
    \smallskip
    Each of these points of view shows that {\it we have associated with $V$ a $G$-torsor over} $k$, i.e. an element $e_V$ of $H^1(k,G)$.
    
    \medskip
    
  \n  7.2. Assume again that char$(k) \neq 2$. As in $\S5$, we may define the virtual quadratic form $q_{4,\varphi}$ and the quadratic forms $q_{5,\varphi}, q_{7,\varphi}, q_{27,\varphi}$; if $\ch(k) \ne 3$, we also have
 $q_{6,\varphi}$.
      Since these forms depend only on  $V$, we may denote them by $q_{4,V}, q_{5,V}, q_{7,V}, q_{27,V}$ and $q_{6,V}$. By (5.5), (5.5) and (5.6), we have 
        
    \medskip
    
    \n {\bf Theorem 7.}
    
    \smallskip
  (7.4)  \quad $q_{27,V} =  \lambda^2q_{5,V} \ + \ \langle1,2\rangle q_{5,V} \ + \ 6  + \langle 2 \rangle $.
    
    \smallskip
    
   (7.5) \quad  $q_{27,V} =    \lambda^2q_{6,V} \ + \ \langle 3 \rangle q_{6,V} \ + \ 6,$  \  \ {\it if}  \ char($k) \neq 3.$

\smallskip

(7.6) \quad $q_{27,V} =  \lambda^2q_{7,V} \ + \ (\langle-1\rangle - \langle 2 \rangle)q_{7,V} \ + \ 7 - \ \langle -2 \rangle $.

    \medskip

\n7.3. {\it  Interpretations of $q_{7,V}$ and} $q_{27,V}$.
      
      \smallskip
      
      \n (a) {\it The case o}f $q_{7,V}$.
      
      \smallskip
 
 Let $V_s$ be the $k_s$-variety deduced from $V$ by the base change $k \to k_s$. The N\'eron-Severi group $L_V$ is equal to the divisor class group $H^1(V_s, \mathcal{O}_{V_s}^\times) $. The map $f \mapsto df/f$ gives a homomorphism $\mathcal{O}_{V_s}^\times \to \Omega^1_{V_s}$; we thus get a homomorphism $L_V  \to H^1(V_s,\Omega^1_{V_s})$, hence also  $k_s \otimes_{\small \Z}L_V \to H^1(V_s,\Omega^1_{V_s})$. Since  $V_s$ is a smooth rational surface, this map is an isomorphism \footnote{ This may be proved by showing that, if it is true for $V$, it is also true for its blow up at any point. 
 One needs to know that the dimension of $H^1(V_s,\Omega^1_{V_s})$ goes up by one, cf. [Ha 77], Chap.V, Exerc.5.3.}. Moreover, it transforms
 the intersection form on $k_s \otimes_{\small \Z}L_V$ into the cup-product: $$ H^1(V_s,\Omega^1_{V_s}) \otimes H^1(V_s,\Omega^1_{V_s}) \to H^2(V_s,\Omega^2_{V_s}) \simeq k_s.$$ By descent, this gives an interpretation of $q_{7,V}$ as the quadratic space $H^1(V,\Omega^1_V)$ endowed
 with its natural cup-product form. (This is almost the same as the Euler characteristic $\chi(V/k)$ in the sense of Levine and Raskit, cf. [LR 18], theorem 3.1; the only difference is that they incorporate the $1$ dimensional quadratic spaces
 $H^0(V,\Omega^0_V) \simeq k $ and  $H^2(V,\Omega^2_V) \simeq k$ in the definition of $\chi(V/k)$.)
 
 \smallskip
 \n (b) {\it The case of} $q_{27,V}$.
 
 \smallskip
 
 That quadratic form is the trace form of the \'etale algebra $A_{27}$ defined by the $\Gamma_k$-set $Y_V$ of the 27 lines. One may also view $A_{27}$
 as the subalgebra of the $G$-Galois algebra $E_V$ of (7.1.1) which is fixed by the subgroup $H \simeq \Weyl(\sf{D}_5)$ defined at the end of \S1.
 
 \medskip
 
 \n 7.4. {\it Questions.}
 
 \smallskip
 
 It would be interesting to be able to compute the quadratic form $q_{7,V}$ (and hence also $q_{27,V}$) directly from the knowledge of a cubic equation $F$ defining $V$; this is probably possible (with computer help) from the interpretation
 of $q_{7,V}$ given in $\S$7.3(a).
 
 A related question is whether the invariants $w_1,...,w_4$ of the
 Appendix satisfy any other identities than those listed at the end of $\S$A.3.
 The difficulty is that not every $\varphi: \Gamma_k \to \sf{E}_6$ comes from a cubic
 surface.
   \vspace{5mm}
    
  \begin{center} {\sc Appendix - Stiefel-Whitney classes}
  \end{center}
  
  The aim of this Appendix is to show how to compute the Stiefel-Whitney classes
  of the linear representations of $G$, and  also those of the quadratic forms $q_{7,\varphi}, q_{27,\varphi},...$ of $\S$5. For instance (cf. $\S$A.4, example 2), we have: $$w_{\rho_{27}} =  \ (1+ew_1+w_2)(1+e^3w_1+w_4)(1+e^6w_2+e^4w_4),$$
 where $e = (-1)$ is the class of 
 $-1$ in $H^1(k)$, and $w_1,...,w_4$ are the basic cohomological invariants of $G$, cf. Theorem 8.
  
  \medskip
  
  We start with recalling some basic definitions (cf. [Ka 84] and [GKT 89]).
  
  \medskip
  
  \n {\bf A.1}. {\it Stiefel-Whitney classes of real  linear representations of $G:$ the group cohomology point of view.}
  
  \smallskip
  
  Here $G$ is any finite group; we denote by $H(G)$ its cohomology
  with coefficients in $\F$$_2$.
  
  Let $\rho:G \to \GL_n(\R)$ be a real linear representation of $G$ of dimension $n$. Let $BG$ and
    $B\GL_n(\R)$ be the classifying spaces of $G$ and $\GL_n(\R)$ in the sense
    of algebraic topology (cf. [Bo 55], $\S$8) - they are well-defined up to homotopy equivalence. The homomorphism $\rho$ induces a map: $BG \to B\GL_n(\R)$,
    hence also a map $\rho^*: H(B\GL_n(\R))$$ \to H(BG)$, where the letter $H$ means cohomology mod 2. Since $\Or_n(\R)$ is a maximal compact subgroup
    of $\GL_n(\R)$, we may identify $H(B\GL_n(\R))$ with $H(B\Or_n(\R))$ which is known to be a polynomial algebra $\F$$_2[w_1,...,w_n]$ in the {\it Stiefel-Whitney
    classes} $w_1,...,w_n$, of degrees $1,...,n$, cf. [Bo 55], $\S$9. Let $w_i(\rho) \in H^i(G)$ be the images of the $w_i$ by $\rho^*$; we put $w_i(\rho) = 1 $ if $i=0$ and $w_i(\rho)
=0$ for $i >n$. The {\it total Stiefel-Whitney class of $\rho $} is
$$  w(\rho) = \sum_{n\geqslant 0} w_i(\rho).$$
These classes only depend on the equivalence class of $\rho$, i.e. of the character $\chi_\rho$ of $\rho$. They enjoy the same properties as the standard Stiefel-Whitney classes.
For instance:

\smallskip

(A1) $ w(\rho_1 \oplus \rho_2) = w(\rho_1)\!\cdot \!w(\rho_2).$

\medskip

 Suppose $n=1$, in which case $\rho(G)$ is contained in $\{\pm1\}$,
i.e. $\rho$ may be identified with an element $[\rho]$ of $H^1(G)= \Hom(G,\Z/2\Z)$. We then have:

\smallskip

(A2) $w(\rho) = 1 + [\rho].$

\smallskip

Note that (A1) and (A2) give a way to compute $w(\rho) $ for every $\rho$
when $G$ is an elementary abelian $2$-group (for instance a  cube), since $\rho$
is then a direct sum of one-dimensional representations.

\smallskip

In dimension $< 4,$ we have:

\smallskip

(A3) $w_1(\rho) = [\det(\rho)],$

\smallskip

(A4) $w_2(\rho)$ is the obstruction to lifting $\rho$ to $\widetilde{\GL_n(\R)}$,
where $\widetilde{\GL_n(\R)}$ is the central extension of $\GL_n(\R)$ by $\{\pm1\}$ characterized by the fact that  a reflection (resp. a product of two commuting distinct reflections) of $\GL_n(\R)$ lifts to an element of order 2 (resp. of order 4) of  $\widetilde{\GL_n(\R)}$.

\smallskip
(A5) $w_3(\rho) = w_1(\rho)\!\cdot\!w_2(\rho) \ + \  \Sq^1w_2(\rho)$,

\smallskip  
\n where $\Sq^1$ is the first Steenrod square operator (i.e. the Bockstein map).

\n {\it Remark.} The definitions above use topology. The reader will find in [GKT 89] a purely algebraic definition of the $w_i(\rho)$, based on functoriality, and properties (A1), (A2), (A4).

\bigskip

   \n {\bf A.2}. {\it Galois Stiefel-Whitney classes.} 
     
  \smallskip
  
  We keep the notation above. 
  
 Let us first recall a general construction. Let $k$ be a field of characteristic $\ne 2$,
  and let $t \in H^1(k,G)$ be a $G$-torsor over $k$. Let $\varphi_t: \Gamma_k \to G$
  be a homomorphism corresponding to $t$. 
  
  Let $z$ be an element of $H(G)$.
  With $t\in H^1(k,G)$ we associate $z_t = \varphi^*_t (z) \in H(k)$. This defines
  a map $H^1(k,G) \to H(k)$, and, when this construction is applied to all the
  extensions of $k$, it gives a {\it cohomological invariant} of $G$ mod 2, i.e. an element of $ \Inv_k(G,\Z/$2$\Z)$, cf. [Se 03], $\S$, Chap.VI. The map
  
 \begin{center} $H(G)  \to \Inv_k(G,\Z/$2$\Z)$ 
 
 \end{center}
 
\n  so defined is one of the main methods for constructing cohomological
  invariants. 
  
  \smallskip

  If we apply this  to the cohomology class $w_i(\rho) \in H^i(G)$, we obtain 
 an element of $\Inv_k^i(G, \Z/$2$\Z)$, which we write $w_{i,\rho}$; one may call it  the {\it $i$-th Galois Stiefel-Whitney invariant}. The total invariant is:
   
 \smallskip

 (A6) \ $w_\rho = \sum_{i \geqslant 0} w_{i,\rho}.$
 
 \medskip
   
  \n {\it Example}. When $G = \sf{S}$$_n$ and $\rho$ is the standard representation
  of degree $n$ of $G$, the $w_{i,\rho}$ are the same as the $\w_i$ of [Se 03], $\S$25.1.
  
  \medskip
  
  \n {\it Remark.} The algebraic properties of the $w_{i,\rho}$ are somewhat simpler
  than those of the $w_i(\rho)$. For instance, property (A5) becomes:
  
  \smallskip
  (A5)$' \ w_{3,\rho} = w_{1,\rho} w_{2,\rho}$,
  
  \smallskip
  
 \n  since Sq$^1$ vanishes on $H^2(k)$. More generally, one has ([Ka 84], Corollary to Theorem 4):
  
  \smallskip
  
$({\rm A}6)' \quad   w_{\rho}=(1+w_{1,\rho})  (1+w_{2,\rho}) (1+w_{4,\rho})  (1+w_{8,\rho})  \cdots $ 
  
    \bigskip
  
  \n {\bf A.3}. {\it The case $G = \Weyl(\sf{E}_6)$. Preliminaries.}
  
  \smallskip
   From now on, we assume that $G = \Weyl(\sf{E}_6)$. Let $\rho_6$ be the natural 6-dimensional representation of $G$, and let $w_i = w_{i,\rho_6} \in \Inv^i_k(G,\Z/$2$\Z)$ be its
    Stiefel-Whitney invariants. We have:

\smallskip

(A7) \ $w_i = 0$ {\it for} $i > 4$.

\smallskip
\n {\it Proof}. Let $C$ be a maximal cube of $G$. The restriction $\rho_6|C$ of $\rho_6$ to $C$ is the direct sum of the standard representation $\rho_4$ of $C$ with a trivial 2-dimensional representation. By (A1), we have

\smallskip

  (A8) \ $w_i(\rho_6|C) = w_i(\rho_4)$ {\it for every} \ $i$,
  
  \smallskip
  
\n  hence $w_i(\rho_6|C) = 0$ for $i > 4$. This shows that the restrictions to $C$ of the $w_i, i > 4$, are $0$. By Theorem 6, this implies (A7).

    \medskip
    \n {\bf Theorem 8.}  {\it The classes $\{w_0, w_1, w_2, w_3, w_4\}$ make up a
    basis of the $H(k)$-module} $\Inv_k(G,\Z/$2$\Z)$.    
    \medskip
    
    \n {\it Proof}.  Let $S$ be the set of the four reflections
    belonging to $C$. By  [Se 03], $\S$16.4, the $H(k)$-algebra $\Inv_k(C,\Z/$2$\Z)$ has a natural basis $a_I$ indexed by the subsets $I$ of $S$. The normalizer $N$ of $C$ acts on $C$ and hence on $\Inv_k(C,\Z/$2$\Z)$; it follows from  (3.10) that it acts transitively
    on the $a_I$ with the same $|I|$. Let us denote by
 $\Inv_k(C,\Z/$2$\Z)^{\rm sym}$ the subalgebra of $\Inv_k(C,\Z/$2$\Z)$ made up of the elements which are fixed under that action.  If $0 \leqslant i \leqslant 4$, let $a_i = \sum_{ |I|=i}a_I$.
 The $a_i$ make up an $H(k)$-basis of   $\Inv_k(C,\Z/$2$\Z)^{\rm sym}$. On the other hand,
 we have:
 
 \smallskip
 (A9) $a_i = w_i|C$.
 
 \smallskip
 
 The restriction map $\Inv$$_k(G,\Z/$2$\Z) \to \Inv$$_k(C,\Z/$2$\Z)$ takes its values in $\Inv_k(C,\Z/$2$\Z)^{\rm sym}$.
 It is injective by Theorem 6, and it is surjective by (A9). 
Hence it is an isomorphism; this proves Theorem 8.
 
 \smallskip
 
 \n {\it Remark.} The proof above is essentially the same as the one
 used in [Se 03], $\S$25.5 to determine $\Inv_k(\sf{S}$$_n)$.
    
    \smallskip
    
    \n {\it Multiplication formulas for the $w_i$.} Before giving these formulas, one more notation is necessary:
    
    If $\alpha$ is an element of $k^\times$, we denote by $(\alpha)_k$, or simply $(\alpha)$, the corresponding element
of $H^1(k) = k^\times/k^{\times 2}$. The case $\alpha = -1$ will be especially
useful; we  write $(-1)_k$ as $e_k$, or simply $e$. We have $x^2=xe$ for  every $x \in H^1(k)$,  or equivalently $x(-x)=0$. More generally, it follows from Milnor's conjecture (now a theorem) that $x^2=xe^d$ for every $x\in H^d(k)$.

\smallskip

 With this notation, we have:
 
 \smallskip $w_i^2 = w_ie^i \ (i\geqslant 0), w_1w_2=w_3, w_1w_3=w_3e, w_1w_4=0, w_2w_4=0, w_3w_4=0.$

    \bigskip
    
    \n {\bf A.4.} {\it Statement of the theorem}.
    
    \smallskip
    
      Let $\rho$ be a real representation of $G$. Our aim is to give
      an explicit formula for the total Stiefel-Whitney class $w_\rho$
      in terms of the character $\chi_\rho$ of $\rho$. Since $w_\rho$ depends only
      on the restriction of $\rho$ to $C$, it will be enough to know the values
      of $\chi_\rho$ on the five classes of involutions of $G$. For $i = 0,...,4$, let
      $g_i$ be an involution of degree $i$, and let $m_i$ be the multiplicity of $ -1$
      as an eigenvalue of $\rho(g_i).$ We have $m_0=0$ and:
      
       \smallskip

      (A10) \ $m_i = \frac{1}{2}(\chi_\rho(g_0)-\chi_\rho(g_i)).$
      
       \medskip

\n {\it Example} 1. If $\rho = \rho_6$, then $m_i=i$ for every $i$.

 \smallskip

\n {\it Example} 2. If $\rho$ is the permutation representation  $\rho_{27}$ of degree 27 given by the action of $G$ on the set $Y$, the $m_i$ are equal to 0, 6, 10, 12, 12. This follows from Lemma 1 of $\S$3.4.

 \smallskip
\n {\it Example} 3. For the representation $\rho_{45}$ of degree 45 given by the action of $G$ on the set of triangles, the $m_i$ are equal to 0, 15, 20, 19, 16. This follows from the decomposition of $\rho_{45}$ in irreducible factors: $1 + \rho_{20} + \rho_{24}$, given in [ATLAS 85], p.26, combined with the values of the characters of $\rho_{20}$ and $\rho_{24}$ given on the next page.

\medskip

\n {\bf Theorem 9.} $w_\rho = 1 + w_1p_1(e)+w_2p_2(e)+w_3p_3(e)+w_4p_4(e)$,

\n{\it where the polynomials $p_i(x) \in \F$$_2[x]$ are characterized by the equations}:

\smallskip
(A11)  \ $x\!\cdot\!p_1(x) = 1+(1+x)^{m_1}$,

(A12) \ $ x^2p_2(x) = 1 + (1+x)^{m_2}$,

(A13)  \ $x^3p_3(x) = 1 + x\!\cdot\!p_1(x)+x^2p_2(x) + (1+x)^{m_3}$,

(A14) \ $x^4p_4(x) = 1 + (1+x)^{m_4}.$

    \medskip
    
  \n   The proof will be given in $\S$A.6.
    
    \medskip
    
   \medskip
  
  \n {\it Remark.} {\it A priori} it is not obvious that the elements $p_i$ of $\F$$_2[x,x^{-1}]$
  defined by the equations (A12), (A13) and (A14) are polynomials. It will be a consequence of the proof. It can also be deduced from the following congruence properties of the $m_i$:

\smallskip

(A15)\ \  $m_2 \equiv 0 $ (mod 2); $m_3 \equiv m_1+m_2$ (mod 4); $m_4 \equiv 2m_2$ (mod 8).

\medskip
\n {\small [Proof of (A15). For $d=2,3,4$, let $C_d$ be a cube of $G$
of order $2^d$. By a well-known property of characters of finite groups,
 the sum $S_d=\sum_{g\in C_d} \chi_\rho(g)$ is divisible by $|C_d| = 2^d$. We have 
 $S_d = \sum_{i=0}^d$ $d\choose i$$z_i,$ where $z_i = \chi_\rho(g_i)$. 
 
 \n This shows that
 $\sum_{i=0}^d$ $d\choose i$$z_i \equiv 0$ (mod $2^d$). By (A10) we have $z_i = z_0-2m_i$. Hence 
  $ \sum_{i=1}^d$ $d\choose i$$m_i \equiv 0$ (mod $2^{d-1}$).
  
 \n For $d=2$, this gives $2m_1+m_2\equiv 0$ (mod 2), i.e., $m_2 \equiv 0 $ (mod 2).
  
 \n For $d=3$, we get $3m_1+3m_2+m_3 \equiv 0$ (mod 4), i.e., $m_3 \equiv m_1+m_2$ (mod 4).

  \n For $d=4$, we get $4m_1+6m_2+4m_3+m_4 \equiv 0$ (mod 8), i.e.,  $8m_1+10m_2+m_4\equiv 0$ (mod 8),
  hence $m_4 \equiv 2m_2$ (mod 8), since $m_2 \equiv 0$ (mod 2).]}
  
\medskip
      
    \n {\bf Corollary}. {\it If $-1$ is a square in $k$, then}:
    
    \medskip
    (A16) \ $w_\rho = (1 + m_1w_1)(1+\frac{m_2}{2}w_2)(1+\frac{m_4}{4}w_4).$
    
    \medskip
    \n {\it Proof.}  Since $e=0$, the $p_i(e)$ are reduced to their constant term.
    Hence:
    
    \medskip
    
  (A16)$'  \ \ w_\rho=  1 + m_1w_1 +$$ m_2\choose2$$ w_2 + ($$m_1\choose3$$+$$m_2\choose3$$ + $$m_3\choose3$$)w_3+$$m_4\choose4$$w_4.$
  
  \medskip
  
 \n   By (A$6)'$,  this gives  $w_\rho = (1+m_1w_1)(1+ $$m_2\choose2$$w_2)(1+$$m_4\choose4$$w_4)$.  The congruences of  (A15) imply that $m_2\choose2$ $\equiv \frac{m_2}{2}$ (mod 2)
  and $m_4\choose4$ $\equiv \frac{m_4}{4}$ (mod 2). Hence (A16). 
  
\bigskip
  \n {\it Example} 1. If $\rho = \rho_6$, one has $p_i(x) = 1$ for $i=1, 2, 3, 4$, and
Theorem 9 gives the basic equation $w_{\rho_6} = 1 + w_1+w_2+w_3+w_4$.
  
  \bigskip
  
  \n {\it Example} 2. If $\rho = \rho_{27}$, one finds:
  
  \smallskip
   \ $p_1(x)= x + x^3 + x^5,$
  
  \smallskip
   \ $p_2(x) = 1+x^6+x^8,$
  
  \smallskip
   \ $p_3(x)= x^3+x^7+x^9,$
  
  \smallskip
   \ $p_4(x)=1+x^4+x^8.$
  
  \smallskip
  
 \n  Hence $w_{i,\rho_{27}}=0$ when $i$ is odd or $i > 12$, and:
   
   \smallskip
      $w_{2,\rho_{27}} = ew_1+w_2$, \quad \ \quad$w_{4,\rho_{27}}\ =e^3w_1+w_4,$
  
  \smallskip
   $w_{6,\rho_{27}}=e^5w_1 + e^3w_3, \quad w_{8,\rho_{27}} \ = e^6w_2+e^4w_4,$
   
   \smallskip
 $ w_{10,\rho_{27}} = e^8w_2+e^7w_3, \ \ w_{12,\rho_{27}} = e^9w_3 + e^8w_4$.  
 
 \medskip
 
 \n In the style of (A$6)'$, this may also be written as:
 
 \smallskip
 
 $w_{\rho_{27}} =  \ (1+ew_1+w_2)(1+e^3w_1+w_4)(1+e^6w_2+e^4w_4).$
 
  \bigskip
  
  \n {\it Example} 3. If $\rho = \rho_{45}$, one finds:
  
  \smallskip
   \ $p_1(x)= \frac{x^{15}-1}{x-1} = 1+ x + ... + x^{14},$
  
  \smallskip
  
   \ $p_2(x) =  x^2 + x^{14}+ x^{18},$
  
  \smallskip 
   \ $p_3(x)= \frac{x^{18}-x^{14}+x^{13}-x^2}{x-1}= (x^2+ ... + x^{12}) + (x^{14}+...+x^{17}),$
  
  \smallskip
  
   \ $p_4(x)=x^{12}.$

  \smallskip
  
 \n  In particular, $w_{20,\rho_{45}} = e^{18}w_2+ e^{17}w_3$ and $w_{i,\rho_{45}} = 0$ when $i>20$.
  
  \bigskip

    \n {\bf A.5.} {\it A preliminary construction $:$ the modified cohomology ring $H(C)'$.}
    
    \smallskip
    
    Let $\rho$ be a real linear representation of $C$. Its Stiefel-Whitney
    classes $w_i(\rho)$ belong to $H(C) \simeq \F_2[x_1,...,x_4]$, cf. $\S$A.1.
    We are going to define them now in the following ring $H(C)'$:
    we introduce a new variable $y$ of degree 1, and we define $H(C)'$ as the quotient
    of $H(C)[y] = \F_2[x_1,...,x_4,y]$ by the ideal generated by the 
    elements  $x_i^2+x_iy$, $i=1,...,4.$ Note that, in $H(C)'$, the  square
    $z^2$ of an element of degree $d > 0$ is given by:
    
      (A17)  \quad $z^2 = zy^d$.
      
    \n   (Proof by induction on $d$, using the fact that $x_i^2=x_iy.)$
  
  \smallskip   
      A down-to-earth description of $H(C)'$ is as follows. For every $I \subset \{1,2,3,4\}$, denote by  $x^I$ the monomial $\prod_{i\in I}x_i$. Then:
      
      \smallskip
      
      (A18) $H(C)'$ {\it is a free  $\F_2[y]$-module with basis the $x^I,$ and multiplication table} 
      
      \smallskip
\n       $x^Ix^J = x^{I\cup J}y^{\mid I\cap J \mid}.$
     
     \smallskip
     
     An interesting feature of $H(C)'$ is that it is {\it universal} for families of
     four elements in $H^1(k)$, where $k$ is any field of characteristic $\ne 2$.
     More precisely:
     
     \smallskip
     
     \n {\bf Lemma 2.} {\it For every family $z_1,...,z_4$ of elements of 
     $H^1(k)$ there exists one and only one homomorphism $\theta: H(C)' \to H(k)$
     such that $\theta(x_i)=z_i$ and $\theta(y) = e$.}
     
     \smallskip
     This follows from the formula $x^2=xe$ for every $x\in H^1(k)$, cf. $\S$A.3.
     
     \
     
     \smallskip
     
     The symmetric group $\sf{S}_4$ acts on $H(C)'$. Its invariants make up a free
     $\F_2[y]$-module $H(C)'^{\rm sym}$ with basis the elementary symmetric functions $s_0,...,s_4$ of the $x_i$:
     
    $s_0 = 1$,
    
    $s_1 = x_1+ x_2+x_3+x_4,$
  
  $s_2 = x_1x_2 + x_1x_3+x_1x_4+x_2x_3+x_2x_4+x_3x_4,$
  
  $s_3 = x_1x_2x_3 + x_1x_2x_4+x_1x_3x_4+x_2x_3x_4,$
  
  $s_4= x_1x_2x_3x_4.$
  
  \smallskip
  
 \n  The $s_i$ are the images in $H(C)'$ of the $w_i|C$ of $\S$A.3, cf. (A9).
 
\medskip

\n {\it Remark}. The natural map $H(C) \to \Inv_k(C)$ of $\S$A.2 can be extended
to $H(C)'$ by requiring that $y$ is mapped to $e_k$. In the special case $k = \R$,
where $H(k) = \F_2[e]$, one finds that $H(C)' \to \Inv_\R(C)$ {\it is an isomorphism, and } $H(C)'^{\rm sym} \simeq \Inv_\R(G)$.

     \bigskip

    \n {\bf A.6.} {\it Proof of Theorem $9$.}
    
    \smallskip
    
    Let $\rho$ be a real linear representation of $C$ whose character 
    is $\sf{S}_4$-invariant, and let $m_1,...,m_4$ be the corresponding integers,
    as defined in $\S$A.4.
    
     Let $w(\rho) \in H(C)$
    be the total Stiefel-Whitney class of $\rho$, cf. $\S$A.1. Let $w(\rho)'$ be the image of $w(\rho)$ in $H(C)'$. Since $w(\rho)'$ is $\sf{S}_4$-invariant, we may write it in a unique way as
    
    \smallskip
 (A19) \quad   $ w(\rho)' = 1 + \sum_{i=1}^4s_ip_i(y), \ with \ \  p_i \in \F$$_2[y],$
 
 \smallskip
 
 \n where the $s_i$ are the elementary symmetric functions of the $x_i$, as above.
 
\medskip

\n {\bf Theorem 10.} {\it The polynomials $p_1, p_2, p_3, p_4$ have the properties}
(A11), (A12), (A13), (A14) \ {\it  of Theorem 9.}
     
     \medskip
 \n {\it Proof.}  Let $C_1$ be a group of order 2, with generator $g$. Its
 cohomology algebra $H(C_1)$ is a polynomial algebra $\F$$_2[x]$. We may enlarge $H(C_1)$ in $H(C_1)'$ as in the previous section, by adding the
 indeterminate $y$, and dividing by the relation $x^2=xy$. The algebra
 so obtained is free of rank 2 over $\F$$_2[y]$ with basis $\{1,x\}$. Note that
 the subalgebra of   $H(C_1)'$ generated by $x$ is isomorphic to $\F$$_2[y]$;
 moreover, if $p$ is a polynomial in one variable, with coefficients in $\F$$_2$,
 we have  
 
  \medskip

 (A20) \ $ x^np(x) = x^np(y)$ \ {\it for every } $n > 0$.
 
  \medskip

   For $i = 1,...,4$, let $f_i: C_1 \to C$ be a homomorphism such that 
   $g_i = f_i(g)$ is an element of $C$ of degree $i$, and let $\rho_i = \rho \circ f_i$.
 The Stiefel-Whitney class $w(\rho_i)'$ is:
 
 \smallskip
 
 (A21) \   $w(\rho_i)'=(1+x)^{m_i}$;
 
 \smallskip
\n  this follows from the identities (A1) and (A2) of $\S$A.1.

\medskip

   Let us now consider separately the four cases $i=1,...,4$.
   
   \medskip

 \n  {\it The case $i=1$.} The map $f_1: C_1 \to C$ defines $f_1^*: H^1(C)' \to H^1(C_1)'$. The fact
 that $g_1=f_1(g)$ is an involution of degree 1 is equivalent to saying that
  $f_1^*$ maps one of the $x_i$ on $x$ and maps the other ones on 0.
 This implies that the images of $s_1,...,s_4$ by $f_1^*$ are respectively
 $x,0,0,0$. Hence the image of $w(\rho)'$ by $f_1^*$ is $1+xp_1(y)$.
 Since that image is $w(\rho_1)'$, we obtain, by (A21):
 
 \smallskip

   $1 + xp_1(y)= (1+x)^{m_1}$,
 
 \smallskip

 \n and by using (A20), we may rewrite this as:
 
 \smallskip

 (A22)  \ $xp_1(x) = 1+ (1+x)^{m_1},$  \ which is the same as (A11).
 
  \medskip

 \n {\it The case $i=2$.}  The homomorphism
 $f_2^*$ maps two of the $x_i$ to $x$ and the other ones to $0$. Hence the
 images of $s_1,...,s_4$ are  $0,x^2,0,0$, and we have:
 
 \smallskip

 (A23) \ $x^2p_2(x) = 1 +(1+x)^{m_2},$ \ which is (A12).
 
 \medskip

 \n {\it The case $i=3$.} Here the images of $s_1,...,s_4$ are $x,x^2,x^3,0$, hence:
 
 \smallskip

 (A24) \  $xp_1(x) + x^2p_2(x)+x^3p_3(x) = 1 + (1+x)^{m_3},$ which is (A13).
 
 \medskip

  \n {\it The case $i=4$.}  Here the images of $s_1,...,s_4$ are $0,0,0,x^4$, hence:
  
  \smallskip

  (A25) \ $x^4p_4(x) = 1+ (1+x)^{m_4}$, \ which is (A14).

 \medskip
\n This concludes the proof of Theorem 10.

\medskip

\n {\it Remark.} The above proof is a kind of {\it interpolation}. It is similar to determining a polynomial through its values at special points. Here
the Stiefel-Whitney class behaves as a fourth degree polynomial with constant term $1$, and we used
its values at the four types of subgroups of order 2.

\bigskip

\n {\it End of the proof of Theorem $9$.}

 \medskip
 
 The natural homomorphism $H(C) \to \Inv_k(C,\Z/$2$\Z)$ can be factored through the homomorphism 
 $H(C)' \to \Inv_k(C,\Z/$2$\Z)$ mapping the element $y$ of $H(C)'$ to $e_k= (-1)$. The images of the $s_i$ are the $w_i$. By applying that homomorphism to $w(\rho)'$, and using Theorem 10, we obtain the formulas of Theorem 9 for the restrictions to $C$ of  $w_\rho, w_1, w_2, w_3, w_4$. By the splitting principle of Theorem 6, this implies Theorem 9.

 \medskip
 \n {\it Remark.} The above proof only uses the following properties of $G$:
 
$ \bullet $ It is a Weyl group in which the maximal cubes have rank $4$
and are conjugate to each other. 

$\bullet$ The involutions of $G$ of the same degree are $G$-conjugate.

\n This shows that {\it Theorem 9 remains true $($with the same polynomials $p_i)$ when $G$ is replaced by the symmetric groups $\sf{S}_8$ or $\sf{S}_9$}. It also remains true for every $\sf{S}$$_n$, the integer $4$ being then replaced
by $[n/2]$; one just needs to extend the list of the $p_i$ beyond $i=4$.

\bigskip
  \n {\bf A.7.} {\it The Stiefel-Whitney classes of the trace forms $q_{T,\varphi}$.}
  
  \smallskip
  
    Let $T$ be a finite $G$-set and let $\varphi: \Gamma_k \to G$ be a continuous
    homomorphism. These data define an \'etale algebra; let $q_{_{T,\varphi}}$ be its trace form, as in Remark 2 of $\S$6.4. Let $w_iq_{_{T,\varphi}}$ be the Stiefel-Whitney classes of $q_{_{T,\varphi}}$ in the sense of Delzant and Milnor (see e.g.
    [Se 03], $\S$17). By a theorem of B. Kahn (cf. [Ka 84] and [Se 03], $\S$25.7),
    these classes are almost the same as those of the linear representation
    $\rho_{_T}: G \to \sf{S}$$_n \to \GL_n(\R)$, where  $n = |T|$, and $G \to \sf{S}$$_n$
    is the homomorphism giving the action of $G$ on $T$. More precisely:
 
 \medskip
 
(A$26) \  \ w_iq_{_{T,\varphi}} = 
\left\{
\begin{array}{ll}
w_{i,\rho_{_T}} \hspace{22mm} {\rm if}\ {\it i} \ {\rm is \ odd} \\

w_{i,\rho_{_T}} + (2)w_{i-1,\rho_{_T}} \ \ {\rm if} \ {\it i} \ {\rm is \ even}.
\end{array}
\right.
$

\medskip

  One may then apply Theorem 9 to $w_{\rho_{_T}}$; the integers $m_i$
  have the following interpretation: if $g_i$ is an involution of degree $i$
  of $G$, $m_i$ is the number of orbits of order 2 of $\{1,g_i\}$ in its action on $T$.
  
  \smallskip
  
  The case of the $G$-set $Y$ with 27 elements is especially simple. We have
  
\smallskip
(A$27) \  wq_{_{Y,\varphi}} = w_{\rho_{27}}$,

\smallskip

\n since all the $w_{i,\rho_{27}}$ are 0 when $i$ is odd. 

The case where $T$ is the
set of the 45 triangles is almost the same. We have

\smallskip

(A$28) \  wq_{_{T,\varphi}} = w_{\rho_{45}} + (2)w_1$.

\smallskip

\n This follows from the fact that the product $(2)e \in H^2(k)$ is 0 since
2 is a sum of two squares. Hence all the terms of $w_{\rho_{45}}$ that are divisible by $e$ disappear after multiplication by $(2)$; the only one that does not is $w_1$.

  \newpage
  \begin{center}  
{\bf Bibliography}
\end{center} 

\bigskip

\n [ATLAS] J.H. Conway, R.T. Curtis, S.P. Norton, R.A. Parker \& R.A. Wilson, {\it Atlas of Finite Groups}, Clarendon Press, Oxford, 1985; second corrected edition, 2003.
     
\n [Bo 55] A. Borel, {\it Topology of Lie groups and characteristic classes}, Bull. A.M.S. 61 (1955), 397-432 (= Coll. Papers, vol. I, no 33).

\n [Bo 68] N. Bourbaki, {\it Groupes et alg\`ebres de Lie, Chap. VI, Syst\`emes de racines}, Hermann, Paris, 1968; English translation, Springer-Verlag, 2002.

\n [Bo 81] \ ----- , {\it Alg\` ebre, Chap. V, Corps commutatifs}, Masson, Paris, 1981;
English translation, Springer-Verlag, 1990.

\n [BS 94] E. Bayer-Fluckiger \& J-P. Serre, {\it Torsions quadratiques et bases normales autoduales}, Amer. J. Math. 116 (1994), 1-63 (= J-P. Serre, Coll. Papers,
vol.~IV, no 163).

\n [Co 88] C.M. Coombes, {\it Every rational surface is separably split}, Comm. Math. Helv. 63 (1988), 305-311.

\n [De 80] M. Demazure, {\it Surfaces de Del Pezzo II-V}, L.N.M. 777 (1980), 23-69.

\n [Do 12] I.V. Dolgachev, {\it Classical Algebraic Geometry $:$ a modern view}, Cambridge Univ. Press, 2012.

\n [GH 19] S. Gille \& C. Hirsch, {\it On the splitting principle for cohomological invariants of reflection groups}, arXiv:1908.08146.

\n [GKT 89] J. Gunarwardena, B. Kahn \& C. Thomas, {\it Stiefel-Whitney classes of real representations of finite groups}, J. Algebra 126 (1989), 327-347.

\n [Ha 77] R. Hartshorne, {\it Algebraic Geometry}, Springer-Verlag, 1977.

\n [Ka 84] B. Kahn, {\it Classes de Stiefel-Whitney de formes quadratiques et de repr\'esentations galoisiennes r\'eelles}, Invent. math. 78 (1984), 223-256.

\n [KMRT 98] M.-A. Knus, A. Merkurjev, M. Rost  \& J-P. Tignol, {\it The Book of Involutions}, AMS Coll. Publ. 44, 1998.

\n [KW 17] J.L. Kass \& K. Wickelgren, {\it An arithmetic count of the lines of a smooth cubic surface}, preprint (2017),  arXiv:1708.01175.

\n [LR 18] M. Levine \& A. Raksit, {\it Motivic Gauss-Bonnet formulas}, arXiv:1808.08385 [math.AG]; to appear in Algebra and Number Theory.

\n [Ma 74] Y. Manin, {\it Cubic Forms, Algebra, Geometry, Arithmetic}, North-Holland, 1974; second edit., 1986.


\n [Sch 58] L. Schl\"afli, {\it An attempt to  determine the twenty-seven lines upon a surface of the third order, and to divide such surfaces into species in reference to the reality of the lines upon the surface}, Quart. J. Pure Appl. Math. 2 (1858), 110-120.

\n [Se 64] J-P. Serre, {\it Cohomologie galoisienne}, LN 5, Springer-Verlag, 1964; fifth revised edit., 1994; English translation, Springer-Verlag, 1997.

\n[Se 03] \ ----- , {\it Cohomological invariants, Witt invariants and trace forms}, notes by Skip Garibaldi, ULECT 28, 1-100, AMS (2003).

\n[Se 18]\  ----- , {\it Cohomological invariants {\rm mod} $2$ of Weyl groups}, Oberwolfach reports 21 (2018), 1284-1286; arXiv:1805.07172.

\n [Sp 52] T.A. Springer, {\it Sur les formes quadratiques d'indice z\'ero}, C.R.A.S. 234 (1952), 1517-1519.

\vspace{25mm}

Eva Bayer-Fluckiger

EPFL-FSB-MATH

Station 8

1015 Lausanne, Switzerland

eva.bayer@epfl.ch

\bigskip

Jean-Pierre Serre

Coll\`ege de France

3 rue d'Ulm

75005 Paris, France

jpserre691@gmail.com

  \end{document}